\documentclass[12pt]{article}

\usepackage{amssymb}
\usepackage{epsfig}
\usepackage{color}
\usepackage[british]{babel}

\newtheorem{theorem}{Theorem}
\newtheorem{lemma}{Lemma}

\newcommand{\rem}{\noindent \textbf{Remark. }}
\newcommand{\rems}{\noindent \textbf{Remarks. }}
\newcommand{\proof}{\noindent \textbf{Proof. }}
\newcommand{\qed}{\square}

\newcommand{\bean}{\begin{eqnarray*}}
\newcommand{\eean}{\end{eqnarray*}}
\newcommand{\bea}{\begin{eqnarray}}
\newcommand{\eea}{\end{eqnarray}}

\newcommand{\inL}{\stackrel{L^1}{\longrightarrow}}
\newcommand{\inLL}{\stackrel{L^2}{\longrightarrow}}
\newcommand{\inLp}{\stackrel{L^p}{\longrightarrow}}
\newcommand{\inLq}{\stackrel{L^q}{\longrightarrow}}

\def\Exp{E}
\def\Pr{P}
\def\R{{\bf R }}
\def\N{{\bf N}}
\def\NN{{\cal N }}
\def\X{{\cal X}}
\def\U{{\cal U}}
\def\M{{\cal M}}
\def\G{{\cal G}}
\def\OO{{\cal O}}
\def\H{{\cal H}}
\def\LL{{\cal L}}
\def\card{{\rm card}}
\def\dist{{\rm dist}}
\def\knng{{k{\rm {\textrm -NNG}}}}
\def\jnng{{j{\rm {\textrm -th~NNG}}}}
\def\onng{{\rm {\textrm ONG}}}
\def\supp{{\rm supp }}
\def\bx{{\bf x}}
\def\by{{\bf y}}
\def\bz{{\bf z}}
\def\bu{{\bf u}}
\def\bv{{\bf v}}
\def\bX{{\bf X}}
\def\bU{{\bf U}}
\def\1{{\bf 1}}
\def\0{{\bf 0}}
\def\po{\preccurlyeq}
\def\potp{\stackrel{\theta,\phi}{\preccurlyeq}}
\def\postar{\preccurlyeq^*}
\def\rd{{\rm d}}
\def\re{{\rm e}}

\begin{document}

\title{Explicit laws of large numbers for random nearest-neighbour type graphs}
\author{Andrew R.~Wade\footnote{e-mail: \texttt{Andrew.Wade@bris.ac.uk}}\\
\normalsize
 Department of Mathematics, University of Bristol,\\
 \normalsize
  University Walk, Bristol BS8 1TW, England.}

\date{February 2007}

\maketitle

\begin{abstract}
Under the unifying umbrella of a general result of Penrose \& Yukich
[\emph{Ann. Appl. Probab.}, (2003) {\bf 13}, 277--303]
we give laws of large numbers (in the $L^p$ sense)
for the total power-weighted length of several nearest-neighbour type graphs
on random point sets
in $\R^d$, $d\in\N$. Some of these results are known; some are new.
We give limiting constants explicitly, where previously they have been evaluated in less
generality or not at all. The graphs we consider include
the $k$-nearest neighbours graph, the Gabriel graph,
the minimal directed spanning forest, and the on-line nearest-neighbour
graph.
\end{abstract}

\vskip 3mm


\vskip 3mm


\section{Introduction}

Graphs constructed
on
 random point sets in $\R^d$
($d \in\N$), formed
 by joining nearby points
according to some deterministic rule, have recently received
considerable interest \cite{p1,st,yu}. Such graphs include the geometric graph,
the minimal spanning tree, and (as studied in this paper)
the nearest-neighbour graph and its relatives. Applications include the modelling
of spatial networks, as well as statistical procedures. 

The graphs in this paper are based on edges between nearest neighbours, sometimes
in some restricted sense. A unifying characteristic of these graphs is {\em stabilization}: roughly speaking,
the configuration of edges around any particular vertex is not affected by changes to the vertex
set outside of some sufficiently large (but finite) ball. Thus these graphs are locally determined in some sense.

A functional of particular interest is the total edge length of the graph, or, more generally, the total
power-weighted edge length (i.e.~the sum of the edge lengths each raised to a given power $\alpha \geq 0$). The 
large-sample asymptotic theory for power-weighted length of
 stabilizing graphs is now well understood;
see e.g.~\cite{kl,p1,p2,py1,py2,st,yu}.

In the present paper we collect several laws of large numbers (LLNs) for total
power-weighted length
 from the family of 
nearest-neighbour type graphs, defined on independent random points
on $\R^d$. We present these results as corollaries to 
a general umbrella theorem of Penrose \& Yukich \cite{py2}. Some of the results 
(for the most common graphs) are known to various extents
 in the literature; others are new. We take a unified
approach which highlights the connections between these results.

In particular, all our results are explicit: we give explicit expressions for limiting constants.
In some cases these constants have been seen
previously in the literature.

Nearest-neighbour graphs and nearest-neighbour
distances in $\R^d$
are of interest in several areas of applied science, including
the social sciences, geography and ecology, where proximity data are often
important (see e.g.~\cite{ko,pi}). Ad-hoc networks, in which nodes scattered in space
are connected according to some geometric rule, are of interest
with respect to various types of communication networks. Quantities of interest
such as overall network throughput may be related to power-weighted length.

In the analysis of multivariate data, in particular via non-parametric
statistics, nearest-neighbour graphs and near-neighbour
distances have found many applications, including
goodness of fit tests, classification, regression, noise estimation,
density estimation, dimension
identification, cluster analysis,
and the two-sample and multi-sample problems;
see for example \cite{bb,bqy,ej,fr,h,hv,t}
and references therein.

\begin{figure}
\begin{center}
\includegraphics[angle=0, width=0.95\textwidth]{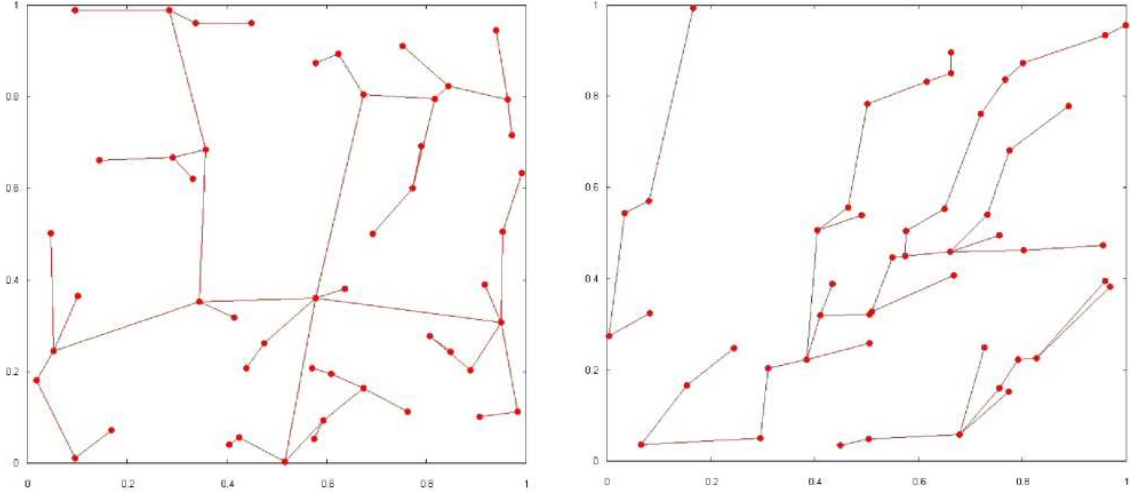}
\end{center}
\caption{Realizations of the ONG (left) and MDSF under $\postar$ (right), each
 on 50 simulated uniform random points in the unit square.}
\label{ongmdst}
\end{figure}

In this paper we give a new LLN for the total power-weighted length
of the {\em on-line nearest-neighbour
graph} (ONG), which is one of the simplest models of network evolution.
We give a detailed description later. In the ONG on a sequence
of points arriving in $\R^d$,
each point after the first is joined by an edge
to its nearest predecessor.
The ONG appeared
in \cite{bbbcr}
as a simple model for the evolution of the Internet graph. Figure \ref{ongmdst} shows a sample
realization of an 
$\onng$.

Recently, graphs with an
`on-line' structure, in which vertices are added
one by one and connected to existing vertices via some rule,
have been the subject of considerable study in relation to the
modelling of real-world networks.  The ONG is one of the simplest 
network evolution models that
captures some of the observed characteristics of real-world networks,
such as spatial structure and sequential growth.

We also consider the {\em minimal directed spanning forest} (MDSF).
The MDSF
is constructed on a partially ordered point set in $\R^d$ by
connecting each point to its nearest neighbour amongst those
points (if any) that precede it in the partial order. If an MDSF is
a tree, it is called a {\em minimal directed spanning tree} (MDST).

The MDST was introduced by Bhatt \& Roy in \cite{br} as a model for drainage 
or communications networks, in $d=2$, with the `coordinatewise' partial order
$\postar$, such that $(x_1,y_1) \postar (x_2,y_2)$ iff $x_1 \leq x_2$ and $y_1 \leq y_2$.
In this version of the MDSF, each point is joined by an edge to its nearest neighbour
in its `south-westerly' quadrant.
In the present paper we give new LLNs for the total power-weighted
length 
for
 a family of MDSFs indexed by partial orderings
 on $\R^2$, which include $\postar$ as a special case.
  Figure \ref{ongmdst}
shows an example of a MDSF under $\postar$. 

\section{Notation and results}

Notions of {\em stabilizing} functionals
of point sets have recently proved to be a useful basis for
establishing limit
theorems
for functionals
of random
point sets in $\R^d$.
In particular, Penrose \& Yukich \cite{py1,py2}
prove general central limit theorems and laws of large numbers
for stabilizing functionals.

The LLNs we give in the present paper
are all derived ultimately from Theorem 2.1 of \cite{py2}, which we restate 
as Theorem
\ref{llnpenyuk} below before we present our results.

In order to describe the result of \cite{py2}, we need to introduce
some notation. Let $d\in \N$.
Let $\| \cdot\|$
be the Euclidean norm on $\R^d$. Write $\card(\X)$ for the cardinality 
of a finite set $\X \subset \R^d$. 
For a locally finite
point set
$\X \subset \R^d$,
$a>0$, and $\by \in
\R^d$, let $\by+a\X$ denote the set $\{ \by + a\bx
: \bx \in \X\}$.
Let $B(\bx;r)$ denote the closed
 Euclidean ball
with centre $\bx \in\R^d$ and radius $r>0$. Let $\0$  denote the origin in $\R^d$.

Let $\xi ( \bx ; \X )$ be a measurable $[0,\infty)$-valued function
defined for all pairs $(\bx , \X)$, where $\X \subset \R^d$ is finite
and $\bx \in \X$. Assume $\xi$ is translation invariant, that is,
 for all $\by \in \R^d$,
 $\xi( \by + \bx; \by+\X) = \xi (\bx ; \X)$.
When $\bx \notin \X$, we abbreviate the notation
$\xi(\bx;\X \cup \{\bx\})$ to $\xi(\bx;\X)$. 
For our applications,  $\xi$ will be {\em homogeneous of order} $\alpha \geq 0$,
that is $\xi ( r \bx ; r \X) = r^\alpha \xi(\bx;\X)$
for all $r>0$, all finite
point sets $\X$, and all $\bx \in \X$.

For any locally finite point set $\X \subset \R^d$ and any $\ell
\in \N$ define 
\bean \xi^+(\X;\ell) := \sup_{k \in \N} \left( \mathrm{ess}
\sup \left\{ \xi(\0;(\X \cap B(\0;\ell)) \cup \mathcal{A}^*) : \mathcal{A}
\in ( \R^d
\backslash B(\0;\ell))^k
\right\} \right) \textrm{, and} \\
 \xi^-(\X;\ell) := \inf_{k \in \N} \left( \mathrm{ess}
\inf \left\{ \xi(\0;(\X \cap B(\0;\ell)) \cup \mathcal{A}^*) : \mathcal{A}
\in ( \R^d
\backslash B(\0;\ell))^k
\right\} \right) , \eean
 where for $\mathcal{A} = (\bx_1,\ldots,\bx_k) \in (\R^d)^k$
we put $\mathcal{A}^*= \{ \bx_1, \ldots, \bx_k \}$ (provided all $k$ vectors are distinct).
 Define the {\em limit} of
$\xi$ on $\X$ by
\[ \xi_{\infty}(\X) := \limsup_{\ell \to \infty}
\xi^+(\X;\ell)  . \] We say the functional $\xi$ \emph{stabilizes} on
$\X$ if
\bean
 \lim_{\ell \to \infty} \xi^+(\X;\ell) = \lim_{\ell \to \infty}
\xi^-(\X;\ell) = \xi_{\infty} (\X)  .
\eean
Stabilization can be interpreted loosely as the property
 that
the value of the functional at a point is unaffected by changes in the
configuration of points at a sufficiently large distance from that point.

Let $f$ be a probability density function on $\R^d$. 
For $n \in \N$
let $\X_n := (\bX_1,\bX_2,\ldots,\bX_n)$ be the point process consisting of $n$ independent
random $d$-vectors 
with common density $f$.
With
probability one, $\X_n$
has
distinct
inter-point distances; hence
 all the
nearest-neighbour type graphs on $\X_n$ that
we consider
are almost
surely unique.

 Let $\H_1$ be a homogeneous Poisson
point process of unit intensity on $\R^d$.
The following general LLN
is due to Penrose \& Yukich, and is obtained from Theorem 2.1 of \cite{py2}
together with equation (2.9) there (the homogeneous case).
\begin{theorem} 
\label{llnpenyuk} Let $q \in \{1,2\}$. Suppose that $\xi$ is homogeneous of order $\alpha$
and
almost surely stabilizes on $\H_1$, with limit
$\xi_{\infty}(\H_1)$.
If $\xi$ satisfies the moments condition \begin{equation}
\label{moms}
 \sup_{n
\in \N} \Exp [ \xi ( n^{1/d} \bX_1 ;n^{1/d}
\X_n  ) ^p  ] < \infty, \end{equation}
 for some $p>q$, then
as $n \to \infty$,
 \bean
n^{-1}
\sum_{\bx \in \X_n} \xi( n^{1/d} \bx; n^{1/d} \X_n )
 \stackrel{L^q}{\longrightarrow}  \Exp [ \xi_\infty ( \H_1)]
\int_{\supp(f)} f(\bx) ^{(d-\alpha)/d} \rd \bx ,
\eean
 and the limit is finite. \end{theorem}

From this result we will derive LLNs for the total power-weighted length
for a collection of nearest-neighbour type graphs.
Let $j \in \N$.
A point $\bx \in \X$ has
a $j$-th {\em nearest neighbour} $\by \in \X \setminus \{ \bx\}$ if
$\card (\{ \bz : \bz \in \X \setminus \{ \bx \}, \|\bz-\bx\| <
\|\by - \bx\| \}) = j-1$.

For all $\bx,\by \in \R^d$ we define the
weight function
\bean
w_\alpha (\bx,\by) := \| \bx-\by \|^\alpha ,\eean
for some fixed parameter $\alpha \geq 0$.
By the total power-weighted edge length of a graph 
with edge set $E$ (where
edges may be directed or undirected),
 we mean the functional
\[ \sum_{ (\bu, \bv) \in E} w_\alpha ( \bu,\bv)
= \sum_{ (\bu,\bv) \in E} \| \bu-\bv\|^\alpha .\]

We will often assume one
of the following conditions on the function $f$ --- either
\begin{itemize}
\item[(C1)] $f$ is supported by a convex polyhedron in $\R^d$ and is bounded away from $0$
and infinity on its support; or
\item[(C2)] for weight exponent $\alpha \in [0,d)$, we require that $\int_{\R^d} f(\bx) ^{(d-\alpha)/d} \rd \bx < \infty$
and $\int_{\R^d} \|\bx\|^r f(\bx) \rd \bx < \infty$ for some $r>d/(d-\alpha)$.
\end{itemize}
In some cases, we take
$f(\bx)=1$ for $\bx \in (0,1)^d$ and $f(\bx)=0$ otherwise, in which case we denote
$\X_n =\U_n=(\bU_1,\bU_2,\ldots,\bU_n)$, the binomial
point process consisting
of $n$ independent uniform random vectors
on $(0,1)^d$.

In the remainder of this section we present our LLNs
derived from Theorem \ref{llnpenyuk}.
Theorems \ref{llnthm}, \ref{nngu}, and \ref{ggthm}
follow directly from Theorem \ref{llnpenyuk} and results in
\cite{py2}, up to evaluation of constants,
while Theorems \ref{onngthm} and \ref{mdstlln} need some more work. 
These results are natural companions, as are their proofs,
which we present in Section \ref{proofs} below; in particular 
the proof of Theorem
\ref{llnthm} is useful for the other proofs.

\subsection{The $k$-nearest neighbours and $j$-th
nearest neighbour graphs}
\label{subsecknng}

Let $j \in \N$. In the
$j$-th nearest-neighbour (directed) graph
on $\X$, denoted by
$\jnng' (\X)$, a directed edge joins each point of $\X$ 
to
its $j$-th nearest-neighbour.

\begin{figure}
\begin{center}
\includegraphics[angle=0, width=0.95\textwidth]{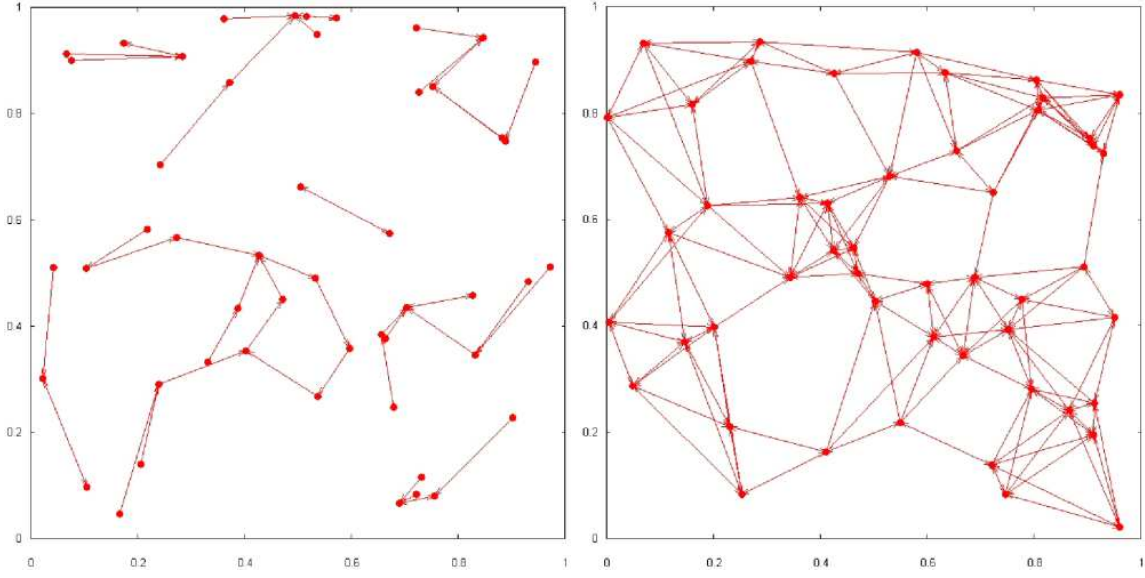}
\end{center}
\caption{Realizations of the  $3$-rd NNG$'$ (left) and
$5$-NNG$'$ (right), each
 on 50 simulated uniform random points in the unit square.}
\label{jkfig}
\end{figure}

Let $k \in \N$. In the $k$-nearest neighbours
(directed) graph on $\X$,
denoted $\knng'(\X)$, a directed edge
joins
each point
of $\X$ 
to each of its first $k$ nearest neighbours in $\X$
(i.e.~each of its $j$-th nearest neighbours for
$j=1,2,\ldots,k$). Clearly the $1$-th NNG$'$ and $1$-NNG$'$
coincide, giving the standard nearest-neighbour (directed) graph. See Figure \ref{jkfig}
for realizations of particular $\jnng'$, $\knng'$.

We also consider the $k$-nearest neighbours
(undirected) graph on $\X$, denoted by
$\knng(\X)$, in which
an undirected edge joins 
$\bx, \by\in \X$ if $\bx$ is one of the first $k$ nearest
neighbours of $\by$,  or $\by$ is one of the first $k$ nearest neighbours
of $\bx$ (or both).

From now on we take the point set $\X$ to be {\em
random}, in particular, for $n \in\N$, we  take $\X=\X_n$.
For $d\in \N$ and $\alpha \geq 0$, let
$\LL^{d,\alpha}_j(\X_n)$, $\LL^{d,\alpha}_{\leq k}(\X_n)$ denote
respectively the total power-weighted edge length of the $j$-th nearest-neighbour (directed)
graph, $k$-nearest neighbours (directed) graph on $\X_n \subset
\R^d$.
Note that
\bea
\label{0818a}
 \LL^{d,\alpha}_{\leq k} (\X_n) = \sum_{j=1}^k \LL^{d,\alpha}_j (\X_n) .\eea

For $d \in \N$, we denote the 
volume of the unit
$d$-ball (see e.g.~(6.50) in \cite{hu}) by
\bea
\label{0818c}
v_d
: = \pi^{d/2} \left[ \Gamma \left( 1+ (d/2) \right) \right]^{-1}.
\eea
 Theorems \ref{llnthm} and \ref{onngthm} below feature
constants $C(d,\alpha,k)$ defined for $d, k \in \N$,  $\alpha \geq 0$ by
\bea
\label{0903c}
C(d,\alpha,k) := v_d^{-\alpha/d} \frac{d}{d+\alpha}
\frac{\Gamma (k+1+(\alpha/d))}{\Gamma(k)}.
\eea

Our first result is Theorem \ref{llnthm} below,
which gives LLNs for
$\LL^{d,\alpha}_j(\X_n)$ and $\LL^{d,\alpha}_{\leq k}(\X_n)$,
with explicit expressions for the limiting
constants; it is the natural starting point for our LLNs for nearest-neighbour
type graphs.
Let $\supp (f)$ denote the support of $f$; under
(C1), $\supp(f)$ is a convex polyhedron, under
(C2) $\supp(f)$ is $\R^d$.

\begin{theorem}
\label{llnthm}
Let $d \in \N$. 
The following results
hold, with 
$p=2$, for $\alpha \geq 0$ if
$f$ satisfies condition (C1), and, with $p=1$, for $\alpha \in [0,d)$ if $f$ satisfies condition
(C2).
\begin{itemize}
\item[(a)]
For $\jnng'$
on $\R^d$ 
we have, as $n \to \infty$,
\bea
\label{0818ff}
 n^{(\alpha-d)/d} \LL^{d,\alpha}_j
(\X_n) \inLp v_d^{-\alpha/d} \frac{\Gamma( j+(\alpha/d))}{\Gamma(j)}
\int_{\supp(f)} f(\bx)^{(d-\alpha)/d} \rd \bx.
\eea
\item[(b)]
For $\knng'$
on $\R^d$
we have, as $n \to \infty$,
\bea
\label{0818g}
 n^{(\alpha-d)/d} \LL^{d,\alpha}_{\leq k}
(\X_n) \inLp C(d,\alpha,k) \int_{\supp(f)} f(\bx)^{(d-\alpha)/d} \rd \bx.
\eea
In particular, as $n\to \infty$,
\bea
\label{0924aa}
 n^{(\alpha-d)/d}  \LL^{d,\alpha}_{\leq k}
(\U_n) \inLp C(d,\alpha,k).
\eea
\end{itemize}
\end{theorem}
\rems
(a) If we use a different norm on $\R^d$ from the Euclidean, Theorem \ref{llnthm}
remains valid with $v_d$ redefined as the volume of the
unit $d$-ball in the chosen norm.

(b) 
Theorem \ref{llnthm}
is essentially contained in Theorem 2.4 of
\cite{py2}, with the  constants evaluated explicitly.
There are several related LLN results in the literature.
Theorem 8.3 of \cite{yu} gives
 LLNs (with complete convergence) for $\LL^{d,1}_{\leq k}(\X_n)$
(see also \cite{mcg}); 
the limiting
constants are not given.
Avram \& Bertsimas (Theorem
7 of \cite{ab}) state a result on the limiting expectation
(and hence the constant in the LLN)
for $\LL^{2,1}_j (\U_n)$, which they attribute to Miles \cite{m} 
(see also p.~101 of \cite{yu}).
The constant in \cite{ab} 
is given as
\[ \frac{1}{2} \pi^{-1/2} \sum_{i=1}^j \frac{\Gamma (i-(1/2))}{\Gamma(i)},\]
which 
simplifies (by induction on $j$)
to $\pi^{-1/2} \Gamma (j+(1/2))/\Gamma(j)$, the $d=2$, $\alpha=1$ case of (\ref{0818ff})
in the case $\X_n =\U_n$.

(c) Related results are the asymptotic 
expectations of
$j$-th nearest neighbour distances in finite
point sets given in \cite{ejs} and \cite{pm}. The
results in \cite{pm} are consistent with the $\alpha=1$
case of our (\ref{0924aa}). The result in \cite{ejs}
includes general $\alpha$ and certain non-uniform densities, although their
conditions on $f$ are more restrictive than our (C1);
the result is consistent with (\ref{0818g}).
Also, \cite{ejs} gives (equation (6.4))
a weak LLN for the empirical
mean $k$-nearest
neighbour {\em distance}.
With Theorem 2.4 of \cite{py2},
the results in \cite{ejs} yield LLNs
 for the total weight of the
$\jnng'$ and $\knng'$ only when $d-1<\alpha<d$
(due to the rates of convergence given in \cite{ejs}). 

(d) Smith \cite{sm} gives, in some sense, expectations of randomly selected edge lengths for nearest-neighbour type
graphs on the homogeneous Poisson point process of unit intensity in $\R^d$, including
the $\jnng'$, nearest-neighbour (undirected) graph, and Gabriel graph. His results coincide with ours
only for the $\jnng'$, since here each vertex contributes a fixed number ($j$) of directed edges: equation (5.4.1) of \cite{sm}
matches the expression for our $C(d,1,k)$. \\

From the results on nearest-neighbour (directed) graphs,
we may obtain results for nearest-neighbour (undirected)
graphs, in which if $\bx$ is a nearest neighbour of $\by$ and vice
versa, then the edge between $\bx$ and $\by$ is counted only once. As an example,
we give the following result.

For  $d \in \N$ and $\alpha \geq 0$ let
$\NN^{d,\alpha}(\X_n)$ denote
the total power-weighted edge length 
 of the nearest-neighbour (undirected)
graph on $\X_n \subset
\R^d$.  
For $d\in \N$, let $\omega_d$
be the volume of the union of two unit $d$-balls with centres unit distance apart in $\R^d$.

\begin{theorem}
\label{nngu}
Suppose that $d \in \N$, $\alpha \geq 0$ and
$f$ satisfies condition (C1). As $n \to \infty$,
\bea
\label{0924x} & &
n^{(\alpha-d)/d} \NN^{d,\alpha} (\X_n) \nonumber\\  
& & \inLL
 \Gamma (1+ (\alpha/d)) \left( v_d^{-\alpha/d}
- \frac{1}{2}  v_d  \omega_d^{-1-(\alpha/d)} \right) \int_{\supp(f)}  f(\bx)^{(d-\alpha)/d} \rd \bx.
\eea
In particular, when $d=2$ we have, for $\alpha \geq 0$
\bea
\label{0924b}
n^{(\alpha-2)/2}  \NN^{2,\alpha} (\U_n) \inLL \Gamma (1+(\alpha/2)) \left( \pi^{-\alpha/2}
- \frac{\pi}{2} \left( \frac{6}{ 8 \pi + 3\sqrt{3}} \right)^{1+(\alpha/2)} \right),
\eea
and when $d=2$, $\alpha=1$, we get
\bea
\label{0913a}
n^{-1/2} \NN^{2,1} (\U_n) \inLL
\frac{1}{2} - \frac{1}{4} \left( \frac{6 \pi}
{8 \pi +3 \sqrt{3}} \right)^{3/2} \approx 0.377508.
\eea
Finally, when $d=1$, $\alpha=1$, we have $\NN^{1,1} (\U_n) \inLL 7/18$ as $n \to \infty$.
\end{theorem}

\rem
A pair of points, each of which is the other's nearest neighbour, is known
as a reciprocal pair. Reciprocal pairs are of interest in ecology (see \cite{pi}).
When $\alpha=0$, $\NN^{d,0}(\X_n)$ counts the number of vertices, minus
one half of the number of reciprocal pairs.
In this case (\ref{0924x})
says 
$n^{-1} \NN^{d,0} (\X_n) \inLL 1 - (v_d/(2 \omega_d))$. This is consistent with results
of Henze \cite{h}
for the fraction of points that are the $\ell$-th
nearest neighbour of their own $k$-th
nearest neighbour; in particular, (see \cite{h} and references therein)
 as $n \to \infty$, the probability that
a point is in a reciprocal pair
tends to $v_d/\omega_d$.

\subsection{The on-line nearest-neighbour graph}
\label{subseconng}

We now consider the {\em on-line nearest-neighbour graph} ($\onng$). 
Let $d \in \N$. Suppose $\bx_1, \bx_2, \ldots$
are points in $(0,1)^d$, arriving
sequentially; for $n\in\N$
form a graph
on vertex set $\{ \bx_1,\ldots,\bx_n\}$
by connecting each
point $\bx_i$, $i=2,3,\ldots,n$ to its nearest neighbour amongst
its predecessors
(i.e.~$\bx_1, \ldots, \bx_{i-1}$),
using the lexicographic ordering on $\R^d$ to break any ties. The
 resulting
tree is
the $\onng$ on $(\bx_1,\bx_2,\ldots,\bx_n)$.

Again, we take our sequence of points to be random.
We restrict our analysis to the case
in which we have independent uniformly distributed points $\bU_1,\bU_2,\ldots$ on $(0,1)^d$.
 For $d \in \N$, $\alpha \geq 0$ and $n \in \N$, let $\OO^{d,\alpha} (\U_n)$
denote the total power-weighted edge length
of the $\onng$ on sequence $\U_n=(\bU_1,\ldots,\bU_n)$.  
The next result gives a new LLN for 
$\OO^{d,\alpha} (\U_n)$
when $\alpha < d$.

\begin{theorem}
\label{onngthm}
Suppose $d \in \N$ and
 $\alpha \in [0,d)$. With $C(d,\alpha,k)$ as given by (\ref{0903c}),
we have that as $n \to \infty$
\bea
\label{0915az}
n^{(\alpha-d)/d} \OO^{d,\alpha} (\U_n)
\inL \frac{d}{d-\alpha} C(d,\alpha,1) = \frac{d}{d-\alpha} v_d^{-\alpha/d} \Gamma (1+(\alpha/d)).
\eea
\end{theorem}

Related results include those on convergence in distribution
of $\OO^{d,\alpha} (\U_n)$,  given in \cite{pw3} for $\alpha >d$ ($\alpha > 1/2$ in the case $d=1$)
and in \cite{p2} in the form of a central limit theorem  for $\alpha \in (0,1/4)$.
Also, the $\onng$ in $d=1$
is related to the `directed linear tree' considered in
\cite{pw2}.

\subsection{The minimal directed spanning forest}
\label{secmdsf}

The {\em minimal directed spanning forest}
 (MDSF) is related to the standard nearest-neighbour (directed)
graph, with the additional constraint that edges can only lie in a given direction. 
In general, the MDSF can be defined as a global optimization problem for directed
graphs on partially ordered sets endowed with a weight function,
and it also admits a local construction; see \cite{br,pw1,pw2}. As above, we consider
the Euclidean setting, where our points lie in $\R^d$.

Suppose that $\X \in \R^d$ is a finite set bearing a partial order
$\preccurlyeq$.
A {\em minimal element}, or {\em sink}, of $\X$ is a vertex
$\bv_0 \in \X$ for which there
exists no $\bv\in \X \setminus \{\bv_0\}$ such that $\bv \preccurlyeq \bv_0$.
Let ${\cal S}$ denote the set of all sinks of $\X$. (Note that ${\cal S}$ cannot be empty.)

For $\bv\in \X$,
we say that $\bu\in \X \setminus \{\bv\}$
 is a \emph{directed nearest neighbour} of $\bv$
if $\bu \preccurlyeq \bv$ and $\| \bv - \bu \| \leq \| \bv - \bu' \|$ for all
$ \bu' \in \X\setminus  \{ \bv \} $ such that $ \bu' \preccurlyeq \bv$.
For each $\bv\in \X \setminus {\cal S}$,
let $n_\bv$ be a directed nearest neighbour of $\bv$ (chosen
arbitrarily if $\bv$ has more than one).
Then (see \cite{pw1}) the directed graph on $\X$
obtained by taking edge set
$
E := \{ (\bv,n_\bv): \bv \in \X \setminus {\cal S} \}
$
is a MDSF of $\X$. Thus, if all edge-weights are distinct, the MDSF
is unique,  and is
 obtained by connecting each non-minimal vertex to its directed
nearest neighbour. In the case where there is a single sink,
the MDSF is a tree (ignoring directedness of edges) and
it is called the {\em minimal directed spanning tree} (MDST).

For what follows, we consider a general type of partial order
on $\R^2$, denoted
$\stackrel{\theta,\phi}{\preccurlyeq}$, specified by the angles
$\theta \in [0 ,2 \pi)$ and $\phi \in (0,\pi ]$.
For $\mathbf{x} \in
\R^2$, let $C_{\theta,\phi}(\mathbf{x})$ be the closed half-cone of angle $\phi$
with
vertex $\mathbf{x}$ and boundaries given by the rays from
$\mathbf{x}$ at angles $\theta$ and $\theta+\phi$, measuring
anticlockwise from the upwards vertical. The partial order is such
that, for $\mathbf{x}_1, \mathbf{x}_2 \in \R^2$,
\bea
 \mathbf{x}_1 \stackrel{\theta,\phi}{\preccurlyeq} \mathbf{x}_2
\textrm{ iff } \mathbf{x}_1 \in C_{\theta,\phi} (\mathbf{x}_2) .
\label{0719} \eea We shall use $\preccurlyeq^*$ as shorthand for
the special case $\stackrel{\pi/2,\pi/2}{\preccurlyeq}$, which is
of particular interest, as in \cite{br}. In this case $(u_1,u_2)
\postar (v_1,v_2)$ 
iff $u_1 \leq
v_1$ and $u_2 \leq v_2$. The symbol $\preccurlyeq$ will denote a
general partial order on $\R^2$. Note that in the case $\phi = \pi$,
(\ref{0719}) does not, in fact, define a partial order
on the whole of $\R^2$, since the antisymmetric property
($\bx \po \by$ and $\by \po \bx$ implies $\bx = \by$) fails; 
however it is, with probability one,
 a true
partial order (in fact, a total order) on the random point sets that we consider. 

We do not permit here the case $\phi=0$,  which
would almost surely give us a disconnected point set.
Nor do
 we allow $\phi \in(\pi, 2\pi]$, since in this case
the
directional relation (\ref{0719})
is not a partial order, since the transitivity property
(if $\bu \preccurlyeq \bv$ and $\bv \preccurlyeq {\bf w}$ then $\bu \preccurlyeq {\bf w}$)
fails for $\phi \in(\pi, 2\pi]$.

Again we take $\X$ to be random;
set $\X= \X_n$, where (as before)
 $\X_n$ is a point
process consisting
of $n$ independent random points on $(0,1)^2$ with common density $f$.
When the partial order
is $\postar$, as in \cite{br}, we also consider the point set $\X_n^0:=\X_n \cup \{ \0\}$ 
(where
$\0$ is the origin in $\R^2$) on which the MDSF is a MDST
rooted at $\0$.

In this random setting, almost surely
each point of $\X$ 
 has  a unique directed nearest neighbour, so
that $\X$ has a unique MDSF.  Denote by $\M^\alpha(\X)$ the total
power-weighted edge length, with weight exponent $\alpha>0$, of the MDSF on $\X$.

Theorem \ref{mdstlln} presents LLNs for $\M^\alpha(\X_n)$
in the uniform case $\X_n=\U_n$.
However, the
proof carries through to other distributions. In particular, if the points
of $\X_n$ are
distributed in
$\R^2$ with a density $f$ that satisfies condition (C1) above,
then (\ref{0728e}) holds with a factor
of $\int_{\supp(f)} f(\bx)^{(2-\alpha)/2} \rd \bx$ introduced into the
right-hand side.

\begin{theorem} \label{mdstlln}
Let $d\in \N$
and
$\alpha \in (0,2)$. Under partial order
$\potp$ with $\theta \in[0,2\pi)$ and $\phi\in(0,\pi]$, 
we have that, as $n \to \infty$,
\bea n^{(\alpha-2)/2}
\M^\alpha ( \U_n ) \inL (2/\phi)^{\alpha/2} \Gamma( 1+(\alpha/2)). \label{0728e} \eea 
Moreover, when the
partial order is $\postar$, (\ref{0728e}) remains true with 
$\U_n$ replaced by $\U_n^0$.
\end{theorem}

\subsection{The Gabriel graph}
\label{gg}

In the {\em Gabriel graph} (see \cite{gs})
on point set $\X \subset \R^d$, two points
are joined by an edge
iff the ball that has the line segment joining those two
points as a diameter contains no other points of $\X$.
The Gabriel
graph has been applied in many of the same contexts as nearest-neighbour graphs;
see for example \cite{t}.

For $d \in \N $ and $\alpha \geq 0$, let
$\G^{d,\alpha}(\X)$ denote
the total power-weighted edge length of the Gabriel 
graph on $\X \subset
\R^d$.
As before, we consider the random point set $\X_n$
with underlying density $f$.
A LLN for $\G^{d,\alpha}(\X_n)$
was given
in \cite{py2}; in the present paper we give the limiting constant explicitly.

\begin{theorem}
\label{ggthm}
Let $d \in \N$ and $\alpha \geq 0$. 
 Suppose that
$f$ satisfies  (C1).
As $n \to \infty$,
\bea
\label{01f}
 n^{(\alpha-d)/d} \G^{d,\alpha}
(\X_n) \inLL v_d^{-\alpha/d} 2^{d+\alpha-1} \Gamma( 1+(\alpha/d)) 
\int_{\supp(f)} f(\bx)^{(d-\alpha)/d} \rd \bx.
\eea
\end{theorem}

\section{Proofs}
\label{proofs}

\subsection{Proof of Theorems \ref{llnthm} and \ref{nngu}} \label{seclln}

For $j \in \N$, let $d_j ( \bx; \X )$ be the (Euclidean) distance from
$\bx$ to its $j$-th nearest
neighbour in $\X \setminus \{\bx\}$,
 if such a neighbour exists, or zero otherwise. 
We will use the following form of
Euler's Gamma integral (see equation 6.1.1 in \cite{as}). For $a,b,c \geq 0$
\bea
\label{0819a}
 \int_0^\infty r^a \re^{-cr^b} \rd r  =
\frac{1}{b} c^{-(a+1)/b
} \Gamma \left( (a+1)/b \right) .
\eea

\noindent
{\bf Proof of Theorem \ref{llnthm}.}
In applying Theorem \ref{llnpenyuk} to the $\jnng'$ and $\knng'$ functionals, we
take $\xi(\bx ; \X_n)$ to be $(d_j(\bx;\X_n))^\alpha$,
 where $\alpha \geq 0$. Then
$\xi$ is translation invariant and homogeneous of order $\alpha$.
It was
shown in Theorem 2.4 of \cite{py2} 
that
 the $\jnng'$ total weight
functional
$\xi$ satisfies
the conditions of Theorem \ref{llnpenyuk} in the following two cases:
(i) with $q=2$, if  $f$ satisfies  (C1), and $\alpha \geq 0$;
and (ii) with $q=1$, if $f$ satisfies  (C2), and $0 \leq \alpha <d$.
(In fact, in \cite{py2} this is proved for the $\knng'$ functional
$\sum_{j=1}^k (d_j(\bx;\X_n))^\alpha$, but this implies
that the conditions also hold for the $\jnng'$ functional
$(d_j(\bx;\X_n))^\alpha$.)

The functional $\xi(\bx;\X_n)
 = (d_j (\bx ; \X_n))^\alpha$
stabilizes on $\H_1$,
with limit $\xi_\infty (\H_1) = (d_j (\0 ; \H_1 ))^\alpha$. Also,
 the moment condition
(\ref{moms}) is satisfied for some $p>1$ (if $f$ satisfies
(C2) and $\alpha<d$) or $p>2$ (if $f$ satisfies (C1)),
and so Theorem \ref{llnpenyuk}, with $q=1$ or $q=2$
respectively, yields (using
the fact that $\xi$ is homogeneous of order $\alpha$)
\bea
 n^{(\alpha/d)-1} \LL^{d,\alpha}_j (\X_n) 
 = n^{-1} \sum_{\bx \in\X_n} n^{\alpha/d} \xi ( \bx ; \X_n )
 =
  n^{-1}
\sum_{\bx \in \X_n} \xi( n^{1/d} \bx; n^{1/d}\X_n) \nonumber\\
 \inLq
\Exp [\xi_\infty(\H_1)] \int_{\supp(f)} f(\bx)^{(d-\alpha)/d} \rd \bx. \label{0728f}
 \eea
We now need to evaluate the expectation on the right-hand side of
(\ref{0728f}). For $r>0$
\bean
\Pr [ \xi_{\infty} ( \H_1 ) >r  ]  & = &
\Pr  [ d_j (\0;\H_1) >r^{1/\alpha}  ]
= \sum_{i=0} ^{j-1}
\Pr  [ \card (\{ B ( \0 ; r^{1/\alpha} ) \cap \H_1 \}) = i
 ]   \\ & = & \sum_{i=0}^{j-1}  \frac{
(v_d r^{d/\alpha})^i}{i!} \exp ( -v_d r^{d/\alpha} ),
\eean
where $v_d$ is given by (\ref{0818c}).
So
\bean
\Exp \left[ \xi_{\infty} (
\H_1 ) \right]  = \int_0^{\infty} \Pr
\left[ \xi_\infty \left( \H_1 \right) > r \right]
 \rd r  = \int_0^{\infty}
\sum_{i=0}^{j-1} \frac{ (v_d
r^{d/\alpha})^i}{i!} \exp ( -v_d r^{d/\alpha}  )  \rd r .
\eean
Interchanging the order of summation and
integration, and using (\ref{0819a}), we obtain
\bea
\label{0819b}
\Exp \left[ \xi_{\infty} (
\H_1 ) \right]
  =   v_d^{-\alpha/d} \frac{\alpha}{d} \sum_{i=0}^{j-1} \frac{\Gamma(i +(\alpha/d))}{\Gamma(i+1)}
= v_d^{-\alpha/d} \frac{\Gamma
\left(j+(\alpha/d) \right)}{\Gamma(j)},
\eea
where the final equality follows by induction on $j$. Then from  (\ref{0818c}),
(\ref{0728f}) and (\ref{0819b}) we obtain the $\jnng'$ result (\ref{0818ff}).
By (\ref{0818a}), the $\knng'$ result (\ref{0818g}) follows
from (\ref{0818ff}) with 
\bean
C(d ,\alpha,k) =
 v_d^{-\alpha/d} \sum_{j=1}^k \frac{\Gamma
\left(j+(\alpha/d) \right)}{\Gamma(j)}
 = v_d^{-\alpha/d} \frac{d}{d+\alpha} \frac{\Gamma(k+1+(\alpha/d))}{\Gamma(k)}. ~\qed
\eean
 
 \noindent
 {\bf Proof of Theorem \ref{nngu}.}
The nearest-neighbour (directed) graph counts
the weights of edges
from points that are nearest neighbours of their
own nearest neighbours twice, while the
nearest-neighbour (undirected) graph
counts such weights only once.

Let $q(\bx;\X)$ be 
the distance from $\bx$
to its nearest neighbour in $\X \setminus \{ \bx\}$
if $\bx$ is a nearest neighbour
of its own nearest neighbour, and zero otherwise. Recall
that $d_1 (\bx;\X)$ is the distance
from $\bx$ to its nearest neighbour in $\X \setminus \{\bx\}$. For $\alpha \geq 0$, 
define
\[ \xi' (\bx;\X) := (d_1 (\bx;\X))^\alpha - \frac{1}{2} ( q(\bx;\X))^\alpha.\]
Then
$\sum_{\bx \in \X} \xi'(\bx,\X)$ is 
the total weight of the nearest-neighbour (undirected) graph on $\X$. Note
that $\xi'$ is translation invariant and homogeneous
of order $\alpha$.

One can check that $\xi'$ is stabilizing on the Poisson process
$\H_1$, using similar arguments to those for the $\jnng'$ and
$\knng'$ functionals. Also (see \cite{py2}) if condition (C1)
holds then $\xi'$ satisfies the moments condition (\ref{moms}) for
some $p>2$, for all $\alpha \geq 0$.

Let ${\bf e}_1$ be a vector
of unit length in $\R^d$. For $d \in \N$,
let $\omega_d:= | B(\0;1) \cup B( {\bf e}_1 ; 1)|$,
the volume of the union of two unit $d$-balls with
centres unit distance apart.

Now we apply Theorem \ref{llnpenyuk}
with $q=2$. We have
\bea
\label{0924c}
 n^{(\alpha/d)-1} \NN^{d,\alpha} (\X_n) =
n^{-1} \sum_{\bx \in \X_n} \xi'(n^{1/d} \bx; n^{1/d} \X_n )
 \nonumber\\
 \inLL  
 \Exp[ \xi'_\infty (\H_1)] \int_{\supp(f)}  f(\bx)^{(d-\alpha)/d} \rd \bx,
\eea
where $\Exp [ \xi'_\infty (\H_1) ] = \Exp[ ( d_1 (\0; \H_1))^\alpha] - (1/2) \Exp [(q (\0; \H_1))^\alpha]$.
Now we need to evaluate $\Exp[(q(\0;\H_1))^\alpha]$. With $\bX$ denoting the nearest
point of $\H_1$ to $\0$,
\bean
 \Pr [ q(\0;\H_1) \in \rd r] & = &
\Pr [ \{|\bX| \in \rd r \} \cap \{ \H_1 \cap ( B(\0;r) \cup B(\bX;r))= \{ \bX \} \} ]
\\
& = & d v_d r^{d-1} \re^{-v_d r^d} \re^{-(\omega_d-v_d)r^d} \rd r
= d v_d r^{d-1} \re^{-\omega_d r^d} \rd r.
\eean
So using (\ref{0819a}) we obtain
 \bea
\label{0924d}
 \Exp [ (q(\0;\H_1))^\alpha] = \int_0^\infty d v_d r^{d-1+\alpha} \re^{-\omega_d r^d} \rd r =
v_d \omega_d^{-1-(\alpha/d)} \Gamma (1+(\alpha/d)).\eea
Then from (\ref{0924c}) with (\ref{0924d}) and the $j=1$ case of
(\ref{0819b})
we obtain (\ref{0924x}). By some calculus, $\omega_2 = (4 \pi/3) +(\sqrt{3}/2)$,
which with the $d=2$ case of (\ref{0924x}) yields (\ref{0924b}); for (\ref{0913a})
note that $\Gamma(3/2)=\pi^{1/2}/2$ (see 6.1.9 in \cite{as}). Finally, we obtain the statement for $\NN^{1,1}(\U_n)$
from the $d=1$ case of (\ref{0924x}) since $\omega_1=3$. $\qed$

\subsection{Proof of Theorem \ref{onngthm}}
\label{onngprf}

In order to obtain our LLN 
 (Theorem \ref{onngthm}
above), we modify the setup of the $\onng$ slightly. Let
$\U_n$ be a {\em marked} random finite point process in $\R^d$,
consisting of $n$ independent uniform random vectors in $(0,1)^d$,
where each point $\bU_i$ of $\U_n$ carries a random mark
$T(\bU_i)$ which is uniformly distributed on $[0,1]$, independent
of the other marks and of the point process $\U_n$. 
Join each point $\bU_i$ of $\U_n$ to its nearest neighbour
amongst those points of $\U_n$
with mark less than $T(\bU_i)$, if there are any such 
points, to obtain a graph that we call the $\onng$ on the
marked point set $\U_n$.
This definition extends to infinite but
locally finite point sets.

Clearly the $\onng$ on the marked point process $\U_n$
has the same distribution as the $\onng$ (with
the first definition) on a sequence $\bU_1,\bU_2,\ldots,\bU_n$
of independent uniform points on $(0,1)^d$.

We apply Theorem \ref{llnpenyuk} to obtain a LLN
for $\OO^{d,\alpha} ( \U_n )$, $\alpha \in [0,d)$. Once again,
the method enables us to evaluate
the limit explicitly.
We
take $f$ 
to be the indicator of $(0,1)^d$.
Define $D(\bx;\X)$ to be
 the distance from point $\bx$ with mark $T(\bx)$ to its nearest
neighbour in $\X$ amongst those points
$\by \in \X$ that have mark $T(\by)$ such that
$T(\by) < T(\bx)$, if such a neighbour exists, or zero otherwise.
We
take $\xi(\bx ; \X)$ to be $(D(\bx;\X))^\alpha$.
Again,
$\xi$ is translation invariant and homogeneous of order $\alpha$.

\begin{lemma}
\label{0919d}
The ONG functional $\xi$ almost surely stabilizes on $\H_1$.
\end{lemma}
\proof
Although the notion of stabilization there is somewhat
different, the same argument as given
at the start of the proof of Theorem 3.6 of \cite{p2} applies. $\qed$

\begin{lemma}
\label{0919a}
Let $d \in \N$, $\alpha \in [0,d)$, and let $p>1$ with $\alpha p<d$.
Then the ONG functional $\xi$ satisfies the moments condition
(\ref{moms}).
\end{lemma}
\proof
Let $T_n$ denote the rank of the mark of $\bU_1$ amongst
the marks of all the points of $\U_n$, so that
$T_n$ is distributed uniformly over
the integers
$1,2,\ldots,n$.
We have, by conditioning
on $T_n$,
\bea
\label{0920a}
\Exp [ (\xi ( n^{1/d} \bU_1 ; n^{1/d} \U_n ))^p ]
& = & n^{-1} \sum_{i=1}^n \Exp [ (d_1 ( n^{1/d} \bU_1; n^{1/d} \U_i )) ^{p\alpha}
] \nonumber\\
& = & n^{-1} \sum_{i=1}^n (n/i)^{p\alpha /d} \Exp [( d_1 (i^{1/d} \bU_1;
i^{1/d} \U_i ))^{p \alpha} ].\eea
It was shown in \cite{pw3} that there exists $C\in(0,\infty)$ such that for all $r>0$
\[ \sup_{i \geq 1} \Pr [ d_1(i^{1/d} \bU_1;
i^{1/d} \U_i ) > r ] \leq C \exp (-r^{1/d}/C).\]
Thus the last expectation in (\ref{0920a}) is bounded by a constant independent of $i$.
So the final expression
in (\ref{0920a}) is bounded by a constant times
\[ n^{(p \alpha -d)/d} \sum_{i=1}^n i^{-p\alpha/d} ,\]
which is uniformly bounded by a constant for $\alpha p <d$. $\qed$ \\

\noindent
{\bf Proof of Theorem \ref{onngthm}.}
Let $d \in \N$. Let $f$ be the indicator of $(0,1)^d$, and
$\xi$ be the $\onng$ functional $\xi(\bx;\U_n)
 = (D (\bx ; \U_n))^\alpha$.
By Lemmas \ref{0919d} and \ref{0919a},
$\xi$
is homogeneous of order $\alpha$,
stabilizing on $\H_1$
with limit $\xi_\infty(\H_1) = ( D( \0 ;\H_1))^\alpha$,
and satisfies the moment condition
(\ref{moms}) for some $p>1$, provided $\alpha < d$. So 
Theorem \ref{llnpenyuk} with $q=1$ implies
\bean
n^{(\alpha/d)-1} \OO^{d,\alpha}(\U_n) =
  n^{-1}
\sum_{\bx \in \U_n} (D (n^{1/d} \bx; n^{1/d}\U_n))^\alpha
\inL
\Exp [\xi_\infty(\H_1)].
 \eean
For $u \in (0,1)$ the points of $\H_1$ with lower mark than $u$ form a
homogeneous Poisson point process of intensity $u$, so by
conditioning on the mark of the point at $\0$,
\bean
\Exp [ \xi_\infty (\H_1)]  =  \int_0^1 \Exp [ (d_1(\0;\H_u))^\alpha ] \rd u 
 =  \int_0^1 u^{-\alpha/d} \Exp [( d_1(\0 ; \H_1))^\alpha ] \rd u = \frac{d}{d-\alpha} C(d,\alpha,1),
\eean
since we saw in the proof of Theorem \ref{llnthm} that $\Exp[ (d_1(\0;\H_1))^\alpha]=C(d,\alpha,1)$. 
$\qed$

\subsection{Proof of Theorem \ref{mdstlln}} \label{seclln2}

In applying Theorem
\ref{llnpenyuk} to the MDSF, we
take 
$f$ to be the indicator
of $(0,1)^2$.
We take $\xi(\mathbf{x} ; \X)$ to be
$(d(\bx;\X))^\alpha$, where
$d(\bx;\X)$ is
 the distance from point $\mathbf{x}$ to its directed nearest
neighbour in $\X \setminus \{\bx\}$, if such a point exists, or
zero otherwise, i.e.
\bea
 \xi(\mathbf{x};\X) = ( d(\mathbf{x};\X) )^\alpha ~~~
{\rm with}~~~
d (\bx ; \X) := \min
\{ \| \bx - \by \| : \by \in \X \setminus \{ \bx \}, \by
\potp \bx \}
\label{0802}
\eea
with the convention that $\min \emptyset = 0$.

We consider the random point set
$\U_n$,
the
binomial point process consisting of $n$ independent uniformly distributed
points on $(0,1)^2$. However, as remarked before the statement of Theorem
\ref{mdstlln}, the result
 (\ref{0728e}) carries through
(with virtually the same proof) to more general point sets $\X_n$.

We need to show that $\xi$ given by (\ref{0802}) satisfies the conditions
of Theorem \ref{llnpenyuk}. As before, $\H_1$ denotes a homogeneous
Poisson process on $\R^2$.

\begin{lemma} \label{stabil} The MDSF functional
$\xi$ given by (\ref{0802}) almost surely stabilizes on
$\mathcal{H}_1$
with limit
$\xi_\infty(\H_1) = (d(\0;\H_1))^\alpha$.
\end{lemma}
\proof
Set $R := d( \0 ; \H_1 )$.
Since $\phi>0$ we have $0< R < \infty$  almost surely.
 But then
for any $\ell >R$, we have $\xi(\0; (\H_1 \cap B(\0;\ell)) \cup
{\cal A}) = R^\alpha$, for any finite ${\cal A} \subset \R^d \setminus
B(\0;\ell)$. Thus $\xi$  stabilizes on $\H_1$ with limit
$\xi_\infty(\H_1) =R^\alpha$.  $\qed$\\

We now give a geometrical lemma.  For $B \subset \R^2$
with $B$ bounded, and for $\bx \in B$,
write $\dist(\bx;\partial B)$ for $\sup\{r: B(\bx;r) \subseteq B\}$,
and for $s >0$, define the region
\bea
\label{Atpdef}
A_{\theta,\phi}(\mathbf{x},s;B) :=
B( \mathbf{x}; s ) \cap B \cap C_{\theta,\phi}(\bx).
\eea
\begin{lemma}
\label{lem0727}
Let $B$ be a  convex bounded set in $\R^2$, and  let $\bx \in B$.
If $A_{\theta,\phi} (\bx,s;B) \cap \partial B(\bx;s) \neq \emptyset$,
and $s > \dist(\bx, \partial B)$, then
$$
|A_{\theta,\phi}(\bx,s;B)| \geq
s
\sin (\phi /2)
 \dist(x,\partial B)  /2.
$$
\end{lemma}
\proof
The condition
 $A_{\theta,\phi} (\bx,s;B) \cap \partial B(\bx;s) \neq \emptyset$
says that there exists $\by \in B \cap C_{\theta,\phi}(\bx)$
with $\|\by - \bx\| = s$.
 The line segment $\bx \by$ is contained in
the cone $C_{\theta,\phi}(\bx)$; take a half-line ${\bf h}$ starting
from $\bx$, at an angle $\phi/2$ to the line segment $\bx \by$
and such that ${\bf h}$ is also contained in $ C_{\theta,\phi}(\bx)$.
Let $\bz$ be the point in $\bf h$ at a distance $\dist(\bx,\partial B)$
from $\bx$. Then the interior of the triangle $\bx \by \bz$ is entirely
contained in  $A_{\theta,\phi}(\bx,s)$, and has area
$s \sin (\phi /2)  \dist(x,\partial B)/2$. $\qed$

\begin{lemma} \label{lem0k715b} Suppose
$\alpha >0$.
Then the MDSF functional $\xi$ given by (\ref{0802}) satisfies the
moments condition (\ref{moms}) for any $p \leq 2 /\alpha$.
\end{lemma}
\proof
Setting $R_n :=(0,n^{1/2})^2$, conditioning on the position of $\bU_1$, we have
\bea
 \Exp  [ \xi  ( n^{1/2}\mathbf{U}_1;n^{1/2} \U_n
 )^p  ]  = n^{-1}
\int_{R_n} \Exp [ ( \xi
(\mathbf{x};n^{1/2} \U_{n-1}) )^p ]  \rd \mathbf{x} 
 .
\label{0728}
\eea
For $\bx \in R_n$ set $m(\bx) := \dist(\bx, \partial R_n)$.
We divide $R_n$ into three regions
\bean
R_n(1) & : = & \{\bx \in R_n: m(\bx) \leq n^{-1/2} \}; ~~~~
R_n(2)  : =  \{\bx \in R_n:  m(\bx) > 1 \};
 \\
R_n(3) &  : = & \{\bx \in R_n: n^{-1/2} < m(\bx) \leq 1 \}.
\eean
 For all $\bx \in R_n$, we have
$\xi(\bx;n^{1/2} \U_{n-1}) \leq (2n)^{\alpha/2}$, and hence,
since $R_n(1)$ has area at most 4, we can
bound the contribution to (\ref{0728}) from $\bx \in R_n(1)$ by
\bea
\label{0728a}
n^{-1} \int_{R_n(1)} \Exp [ ( \xi
(\bx;n^{1/2}\U_{n-1}))^p ] \rd \bx
\leq  4 n^{-1} (2n)^{p\alpha /2} = 2^{2+ p\alpha/2}  n^{(p \alpha -2)/2},
\eea
which is bounded if $p \alpha \leq 2$.
Now, for $\mathbf{x} \in R_n$, with $A_{\theta,\phi}$ defined
at (\ref{Atpdef}),
we have
\bea
\Pr [  d
(\bx; n^{1/2} \U_{n-1})  > s ] & \leq &
\Pr [ n^{1/2} \U_{n-1}
 \cap A_{\theta,\phi}(\bx,s;R_n) = \emptyset ]
\nonumber \\
 & = &
\left(1 - \frac{|A_{\theta,\phi}(\bx,s;R_n)| }{n} \right)^{n-1}
\nonumber \\
& \leq &  \exp( 1 - |A_{\theta,\phi}(\bx,s;R_n)| ),
\label{0728b} \eea
since $|A_{\theta,\phi}(\bx,s;R_n)|\leq n$.
For $\bx \in R_n$ and $s>m(\bx)$, by Lemma \ref{lem0727} we have
 $$
| A_{\theta,\phi}(\mathbf{x},s;R_n)
 | \geq s \sin(\phi/2 )  m(\bx)/2 ~~~{\rm if} ~~
A_{\theta,\phi} (\bx,s;R_n) \cap \partial B(\bx;s)
\neq \emptyset,
$$
and also
$$
\Pr [ d(\bx;n^{1/2} \U_{n-1}) > s ] = 0 ~~~{\rm if}~~~
A_{\theta,\phi} (\bx,s;R_n) \cap \partial B(\bx;s)
= \emptyset.
$$
 For $s \leq m(\bx)$,
we have that
$
 |
A_{\theta,\phi}(\mathbf{x},s;R_n) | = s^2
(\phi/2)  \geq s^2 \sin (\phi/2).
$
Combining these observations and (\ref{0728b}),
 we obtain for all $\bx \in R_n$  and $s >0$ that
\bean
\Pr  [  d
(\bx; n^{1/2} \U_{n-1})  > s ] & \leq & 
\exp  ( 1 -  (s/2) \min (s, m(\bx)) \sin (\phi/2) 
), ~~~ \bx \in R_n.
\eean
Setting $c =(1/2) \sin(\phi/2)$, we therefore have for $\mathbf{x} \in R_n$
that 
\bea
  \Exp [  (\xi (\mathbf{x};n^{1/2} \U_{n-1}))^p
]
 = \int_0^\infty \Pr  [ d(\bx; n^{1/2} \U_{n-1})
> r^{1/(\alpha p)}  ] \rd r
\nonumber\\
 \leq 
  \int_0^{m(\bx)^{\alpha p}} 
 \exp {   (1 - c r^{2/(\alpha p)}  ) } \rd r
 +
   \int_{m(\bx)^{\alpha p }}^\infty 
\exp {  (1 - c m(\bx)r^{1/(\alpha p)}  ) } \rd r
\nonumber \\
 =  O(1) + \alpha p m(\bx)^{-p\alpha} \int_{m(\bx)^2}^\infty
\re^{1-cu} \alpha p u^{\alpha p-1} \rd u
 =  O(1) + O(m(\bx)^{-\alpha p}).
\label{0728d}
\eea
For $\bx \in R_n(2)$, this bound is $O(1)$, and the area
of $R_n(2)$ is less than $n$, so that the contribution
to (\ref{0728}) from $R_n(2)$ satisfies
\bea
\label{0728c}
\limsup_{n \to \infty} n^{-1}
\int_{R_n(2)}   \Exp [ ( \xi
(\mathbf{x};n^{1/2} \U_{n-1}) )^p ]
\rd \mathbf{x}
< \infty.
\eea
Finally, by (\ref{0728d}), there is a constant $C \in (0,\infty)$ such
that 
the contribution to (\ref{0728}) from $R_n(3)$ satisfies
\bean
n^{-1} \int_{R_n(3)} \Exp [ ( \xi
(\mathbf{x};n^{1/2} \U_{n-1}) )^p ] \rd \mathbf{x}
&\leq &
  C n^{-1/2} \int_{y=n^{-1/2}}^1 y^{-\alpha p} \rd y
\\
& \leq & C n^{-1/2} \max \{ \log n , n^{(\alpha p -1)/2} \},
\eean
which is bounded if $\alpha p \leq 2$.
Combined with the bounds in (\ref{0728a}) and (\ref{0728c}),
this shows that the expression (\ref{0728})
is uniformly bounded, provided $\alpha p \leq 2$.
$\qed$\\

For $k \in \N$, and for $a<b$ and $c<d$
 let $\U_{k,(a,b] \times (c,d]}$ denote the point process
consisting of
$k$ independent  random vectors uniformly distributed
on the rectangle $(a,b] \times (c,d]$.
Before proceeding further,
we recall
that
if $M(\X)$ denotes the number of minimal elements, under partial
order
$\postar$, of a point
set $\X \subset \R^2 $, then
\bea
\Exp [ M(\U_{k,(a,b]\times (c,d]}) ] =
\Exp [ M(\U_k) ]
 = 1 + (1/2) + \cdots + (1/k) \leq 1 + \log k.
\label{harmonicbd}
\eea
The first equality in (\ref{harmonicbd})  comes
from some obvious scaling which shows that the
distribution of
$ M(\U_{k,(a,b]\times (c,d]}) $ does not depend on $a,b,c,d$.
For  the second equality
in (\ref{harmonicbd}),
see e.g.~\cite{bns}.\\

\noindent
{\bf Proof of Theorem \ref{mdstlln}.}
Suppose $\alpha \in (0,2)$,
and set $f$ to be the indicator of 
$(0,1)^2$.
 By Lemmas
\ref{stabil} and \ref{lem0k715b} the functional $\xi$,
given at (\ref{0802}), satisfies the
conditions of Theorem \ref{llnpenyuk} with $p= 2/\alpha$
and $q=1$.
 So by Theorem \ref{llnpenyuk},
we have 
\bea
n^{(\alpha/2)-1} \M^\alpha(\U_n) =
  n^{-1}
\sum_{\bx \in \U_n} \xi (n^{1/2}\mathbf{x};n^{1/2}\U_n)
\inL
\Exp [ \xi_\infty(\H_1)].
\label{0728fm}
 \eea
Since the disk sector $C_{\theta,\phi}(\mathbf{x}) \cap B(\bx;r)$
has area $(\phi/2) r^2$, by Lemma \ref{stabil} we have
 \bean \Pr [ \xi_{\infty} (
\H_1 ) >s  ] & = & \Pr  [ \H_1 \cap
C_{\theta,\phi}(\mathbf{0}) \cap B(\0;s^{1/\alpha}) = \emptyset  ] =
\exp (-(\phi/2) s^{2/\alpha} ).
\eean
Hence the limit in (\ref{0728fm}) is, using (\ref{0819a}),
 \[
  \Exp \left[ \xi_{\infty} ( \H_1 ) \right]  =
\int_0^{\infty} \Pr \left[ \xi_\infty \left( \H_1 \right)
> s \right] \rd s
= \alpha 2^{(\alpha-2)/2} \phi^{-\alpha/2} \Gamma( \alpha/2 ),
\]
and this gives us (\ref{0728e}). Finally, in the case where $\potp$ is $\postar$,
(\ref{0728e}) remains true when $\U_n$ is replaced by $\U_n^0$, since
\bea
\label{0806a}
 \Exp [n^{(\alpha/2)-1} | \M^\alpha (\U_n^0) - \M^\alpha (\U_n)| ]
\leq 2^{\alpha/2} n^{(\alpha/2)-1} \Exp [ M(\U_n)] ,
\eea
where $M(\U_n)$ denotes the number of $\postar$-minimal elements of $\U_n$. By (\ref{harmonicbd}),
$\Exp[ M(\U_n)] \leq 1+\log n$, and hence the right-hand side of (\ref{0806a})
tends to 0 as $n \to \infty$ for $\alpha<2$. 
$\qed$

\subsection{Proof of Theorem \ref{ggthm}} \label{secgg}

\noindent
{\bf Proof of Theorem \ref{ggthm}.}
In applying Theorem \ref{llnpenyuk} to the Gabriel graph, we
take $\xi(\bx ; \X_n)$ to be 
{\em one half} of the total $\alpha$ power-weighted
length
of all the edges incident to $\bx$ in the Gabriel graph on $\X_n \cup \{\bx\}$; the factor
of one half prevents double counting.
As stated in \cite{py2} (Section 2.3(e)),
$\xi$ is translation invariant, homogeneous of order $\alpha$ and
 stabilizing on $\H_1$, and 
 if the function $f$ satisfies condition (C1) 
then the moment condition
(\ref{moms}) is satisfied for some $p>2$.
 So by Theorem \ref{llnpenyuk} with $q=2$, 
\bea
 n^{(\alpha/d)-1} \G^{d,\alpha} (\X_n) =
  n^{-1}
\sum_{\bx \in \X_n} \xi( n^{1/d} \bx; n^{1/d}\X_n) \nonumber\\
 \inLL
\Exp [\xi_\infty(\H_1)] \int_{\supp(f)} f(\bx)^{(d-\alpha)/d} \rd \bx. \label{08f}
 \eea
We need to evaluate the expectation on the right-hand side of
(\ref{08f}). The net contribution
from a vertex at $\0$ to the total weight of the Gabriel graph on $\H_1$ is 
\bea
\label{01a}
 \frac{1}{2} \sum_{k=1}^\infty (d_k (\0;\H_1))^\alpha \cdot \1_{E_k},\eea
where the factor of one half ensures that edges are not counted twice,
$d_k(\0;\H_1)$ is the distance from $\0$ to its $k$-th nearest neighbour
in $\H_1$, and $E_k$ denotes the event that $\0$ and its $k$-th nearest neighbour
in $\H_1$ are joined by an edge in the Gabriel graph.

Given that the point $\bx \in \H_1$ is the $k$-th nearest neighbour
of $\0$, an edge between $\bx$ and $\0$ exists in the Gabriel graph iff
the ball with $\0$ and $\bx$ diametrically opposed contains none of the other
$k-1$ points of $\H_1$ that are uniformly distributed in the ball centre $\0$ and
radius $\|\bx\|$. Thus for $k \in \N$,
\bea
\label{01b}
 \Pr[E_k] = \left( \frac{v_d r^d - v_d(r/2)^d}{v_d r^d} \right)^{k-1}
= \left( 1 - 2^{-d} \right)^{k-1}.\eea
So from (\ref{01a}) and (\ref{01b}) we have
\bean \Exp[\xi_\infty(\H_1)] & = & \frac{1}{2} \sum_{k=1}^\infty 
\left( 1 - 2^{-d} \right)^{k-1} \Exp[ (d_k (\0;\H_1))^\alpha] \\
& = & \frac{1}{2} \sum_{k=1}^\infty 
\left( 1 - 2^{-d} \right)^{k-1} v_d^{-\alpha/d} \frac{\Gamma(k+(\alpha/d))}{\Gamma(k)},
\eean
by (\ref{0819b}). But by properties of Gauss hypergeometric series (see 15.1.1 and 15.1.8 of \cite{as})
\[ \sum_{k=1}^\infty 
\left( 1 - 2^{-d} \right)^{k-1} \frac{\Gamma(k+(\alpha/d))}{\Gamma(k)}
= \Gamma (1+(\alpha/d)) 2^{d+\alpha}.\]
Then with (\ref{08f}) the proof is complete. $\qed$

 \begin{center}
 {\bf Acknowledgements}
 \end{center}
  
Some of this work was done when AW was
at the University of Durham, supported by an EPSRC
doctoral training account, and at the University of Bath. AW thanks Mathew Penrose for very helpful
discussions and comments, and two anonymous referees, whose comments have led to an improved presentation.

\end{document}